\documentclass[12pt]{amsart}
\usepackage{amsmath,amsthm,amsfonts,amscd,amssymb,pb-diagram}

\theoremstyle{plain}
\newtheorem{Thm}{Theorem}[section]
\newtheorem{Lem}[Thm]{Lemma}

\newtheorem{Cor}[Thm]{Corollary}

\newtheorem{Def/Prop}[Thm]{Defnition/Proposition}
\newtheorem{Rem}{Remark}

\theoremstyle{definition}
\newtheorem{Def}[Thm]{Definition}

\def\ra{\rightarrow}

\def\P{\mathbb{P}}

\def\O{\mathcal{O}}

\def\et{\acute{e}t}

\def\til{\widetilde}

\author{Daniel Reuben Krashen}

\title{Moduli of \'etale subalgebras in an Azumaya algebra}

\begin{document}

\begin{abstract}
Let $A$ be a sheaf of Azumaya algebras over a Noetherian base $S$. In
this paper we describe using generalized Severi-Brauer varieties, a
quasi-projective moduli space parametrizing sheaves of \`etale
subalgebras of $A$.

In the case that $S$ is the spectrum of a field, we study the geometry
of this moduli space, and show that in certain cases it is a rational
variety.
\end{abstract}

\maketitle

\section{Preliminaries}

Let $S$ be a Noetherian scheme, and let $A$ be a sheaf of Azumaya
algebras over $S$. Our goal is to study the functor $\et(A)$, which
associates to every $S$-scheme $X$, the set of sheaves of commutative
\'etale subalgebras of $A_X$. We will show that this functor is
representable by a scheme which may be described in terms of the
generalized Severi-Brauer variety of $A$.

Unless said otherwise, all products are fiber products over $S$. If
$X$ is an $S$-scheme with structure morphism $f : X \ra S$, then we
write $A_X$ for the sheaf of $\O_X$-algebras $f^*(A)$. For a
$S$-scheme $Y$, we occasionally write $Y_X$ for $Y \times X$, thought
of as an $X$-scheme.

Every sheaf of \'etale subalgebras may be assigned a discrete
invariant, which we call its type, and therefore our moduli scheme is
acually a disjoint union of other moduli spaces.

To begin, let us define the notion of type. 

\begin{Def}
Let $R$ be a local ring, and $B/R$ an Azumaya algebra. If $e \in B$ is
an idempotent, we define the rank of $e$, denoted $r(e)$ to be the
reduced rank of the right ideal $eB$. Recall that the reduced rank of
a $B$-module $M$ is the rank of $M$ as an $R$-module divided by the
degree of $B$ (\cite{BofInv}).
\end{Def}

Let $E$ be a sheaf of \'etale subalgebras of $A/S$, and let $p \in
S$. Let $R$ be the local ring of $p$ in the \'etale topology (so that
$R$ is a strictly Henselian local ring). Then taking \'etale stalks,
we see that $E_p$ is an \'etale subalgebra of $A_p/R$, and it follows
that
$$E_p = \oplus_{i=1}^k R e_i,$$ for a uniquely defined collection of
idempotents $e_i$, which are each minimal idempotents in $S_p$.

\begin{Def}
The type of $E$ at the point $p$ is the unordered collection of
positive integers $[r(e_1), \ldots, r(e_m)]$.
\end{Def}

\begin{Def}
We say that $E$ has type $[n_1, \ldots, n_m]$ if if has this type for
each point $p \in S$.
\end{Def}

\begin{Rem} \label{idempotents}
Since $1 = \sum e_i$, the ideals $I_i = e_i A$ span $A$. Further it is
easy to see that the ideals $I_i$ are linearly independent since $e_ia
= e_jb$ implies $e_i a = e_i e_i a = e_i e_j b = 0$. We therefore know
that the numbers making up the type of $E$ give a partition of
$deg(A_p)$.
\end{Rem}

Some additional notation for partitions will be useful. Let $\rho =
[n_1, \ldots, n_m]$. For a positive integer $i$, let $\rho(i)$ be the
number of occurences of $i$ in $\rho$. Let $S(\rho)$ be the set of
distinct integers $n_i$ occuring in $\rho$, and let $N(\rho) =
|S(\rho)|$. Let
$$\ell(\rho) = \underset{i \in S(\rho)}{\sum} \rho(i) = m$$
be the length of the partition.

\section{Moduli spaces of \'etale subalgebras}

Suppose $A/S$ is an sheaf of Azumaya algebras, and suppose $S$ is a
connected, Noetherian scheme.  Let $\rho = [n_1, \ldots, n_m]$ be a
partition of $n = deg(A)$. Let $\et_{\rho}(A)$ be the functor which
associates to every $S$ scheme $X$ the set of \'etale subalgebras of
$A_X$ of type $\rho$. That is, if $X$ has structure map $f: X \ra S$, 
$$\et_{\rho}(A)(X) = 
\left\{
\begin{matrix}
\text{sub-$\O_X$-modules } \\
E \subset f^*A 
\end{matrix}
\left|
\begin{matrix}
\text{$E$ is a sheaf of commutative \'etale} \\
\text{subalgebras of $f^*A$ of type $\rho$}
\end{matrix}
\right. \right\}$$

Our first goal will be to describe the scheme which represents this
functor. We use the following notation: 
$$V_\rho(A) = \underset{i \in S(\rho)}{\prod} V_i(A)^{\rho(i)},$$
where $V_i(A)$ is the $i$'th generalized Severi-Brauer variety of $A$
(\cite{Blanchet}), which parametrizes right ideals of $A$ which are
locally direct summands of reduced rank $i$.

We define $V_\rho(A)^\circ$ to be the open subscheme parametrizing
ideals which are linearly independent. That is to say, for a
$S$-scheme $X$, if $I_1, \ldots, I_{\ell{\rho}}$ is a collection of sheaves of
ideals in $A_X$, representing a point in $V_\rho(A)(X)$, then by
definition, this point lies in $V_\rho(A)^\circ$ if and only if
$\oplus I_i = A$.

Let $S_\rho$ be the subgroup $\prod_{i \in S(\rho)} S_{\rho(i)}$ of
the symmetric group $S_n$. For each $i$, we have an action of
$S_{\rho(i)}$ on $V_i(A)^{\rho(i)}$ by permuting the factors. This
induces an action of $S_{\rho}$ on $V_{\rho}(A)$, and on
$V_{\rho}(A)^\circ$. Denote the quotients of these actions by
$S^{\rho}V(A)$ and $S^{\rho}V(A)^{\circ}$ respectively. We note that
since the action on $V_{\rho}(A)^{\circ}$ is free, the quotient morphism
$$V_{\rho}(A)^{\circ} \ra S^{\rho}V(A)^{\circ}$$ is an \'etale
morphism which is a Galois covering with group $S_\rho$.

\begin{Thm} \label{moduli}
Let $\rho = [n_1, \ldots, n_m]$ be a partition of $n$. Then the
functor $\et(A)_\rho$ is represented by the scheme $S^{\rho}V(A)$.
\end{Thm}
\begin{proof}
To begin, we first note that both $\et(A)_\rho$ and the functor
represented by $S^{\rho}V(A)$ are sheaves in the \'etale
topology. Therefore, to show that these functors are naturally
isomorphic, it suffices to construct a natural transformation $\psi :
S^{\rho}V(A)^{\circ} \ra \et(A)_\rho$, and then show that this
morphism induces isomorphisms on the level of stalks.

Let $X$ be an $S$-scheme, and let $p : X \ra S^{\rho}V(A)^{\circ}$. To
define $\psi(X)(p)$, since both functors are \'etale sheaves, it
suffices to define it on an \'etale cover of $X$. Let $\til{X}$ be the
pullback in the diagram
\begin{equation} \label{pullback}
\begin{diagram}
\node{\til{X}} \arrow{e} \arrow{s}
  \node{V_{\rho}(A)^{\circ}} \arrow{s,r}{\pi} \\
\node{X} \arrow{e,t}{p} \node{S^{\rho}V(A)^{\circ}}
\end{diagram}
\end{equation}

Since the quotient morphism $\pi$ is \'etale, so is the morphism
$\til{X} \ra X$. Therefore we see that after passing to an \'etale
cover, and replacing $X$ by $\til{X}$, we may assume that $p = \pi(q)$
for some $q \in V_{\rho}A^{\circ}(X)$. Passing to another cover, we
may also assume that $X = Spec(R)$.

Since $p = \pi(q)$, we may find right ideals $I_1, \ldots,
I_{\ell(\rho)}$ of $A_R$ such that $\oplus I_i = A_R$, which represent
$q$. Writing
$$1 = \sum e_i, \ \ e_i \in I_i,$$ we define $E_p = \oplus e_i
R$. This is a split \'etale extension of $R$, which is a subalgebra of
$A$, and we set $\psi(p) = E_p$. One may check that this defines a
morphism of sheaves. Note that this definition with respect to an
\'etale cover gives a general definition since the association $(I_1,
\ldots, I_{\ell}) \mapsto E_p$ is $S_\rho$ invariant.

To see that $\psi$ is an isomorphism, it suffices to check that it is
an isomorphism on \'etale stalks. In other words, we may restrict to
the case that $X = Spec(R)$, where $R$ is a strictly Henselian local
ring. 

We first show that $\psi$ is injective. Suppose $E$ is an \'etale
subalgebra of $A_R$ of type $\rho$. Since $R$ is strictly Henselian,
we have 
$$E = \underset{i \in S(\rho)}{\oplus} \underset{j =
1}{\overset{\rho(i)}{\oplus}} e_{i,j} R.$$ By definition, since the
type of $E$ is $\rho$, if we let $I_{i,j} = e_{i,j}A$, then we the
tuple of ideals $(I_{i,j})$ defines a point $q \in
V_\rho(A)^{\circ}(R)$. Further, since $\sum e_{i,j} = 1$, we actually have $q
\in V_\rho(A)^{\circ}(R)$. If we let $p = \pi(q)$, then tracing
through the above map yeilds $\psi(R)(p) = E$. Therfore $\psi$ is
surjective.

To see that it is injective, we suppose that we have points $p, p' \in
S^{\rho}V(A)^{\circ}(R)$. By forming the pullbacks as in equation
\ref{pullback}, since $R$ is strictly Henselian, we immediately find
that in each case, because $\til{X}$ is an \'etale cover of $X$, it is
a split \'etale extension, and hence we have sections. This means we may write
$$p = \pi(I_1, \ldots, I_{\ell(\rho)}), \ \ p' = \pi(I_1', \ldots,
I_{\ell(\rho)}').$$ Note that in order to show that $p = p'$ is
suffices to prove that the ideals are equal after reordering. Now, if
$E_p = E_{p'}$, then both rings have the same minimal
idempotents. However, by remark \ref{idempotents}, the ideals are
generated by these idempotents. Therefore, the ideals concide after
reordering, and we are done.
\end{proof}

Since we now know that the functor $\et_\rho(A)$ is representable, we
will abuse notation slightly and refer to it and the representing
variety by the same name.

\begin{Def}
$\et(A)$ is the disjoint union of the schemes $\et_{\rho}(A)$ as
$\rho$ ranges over all the partitions of $n = deg(A)$.
\end{Def}

\begin{Cor}
The functor which associates to any $S$-scheme $X$ the set of \'etale
subalgebras of $A_X$ is representable by $\et(A)$.
\end{Cor}

\begin{Rem} \label{grass}
By associating to an \'etale subalgebra $E \subset A_X$ its underlying
module, we obtain a natural transformation to the Grassmannian
functor. 
\end{Rem}

\section{Subfields of central simple algebras}

In this section we specialize to the case where $S = Spec(k)$ for a
field $k$, and $A$ is a central simple $k$-algebra. We use the
notation in this section that tensor products are always taken over
$k$, unless specified otherwise. If $E$ is an \'etale subalgebra of
$A$, then taking the \'etale stalk at $Spec(k)$ amounts to extending
scalars to the seperable closure $k^{sep}$ of $k$. Let $G$ be the
absolute Galois group of $k^{sep}$ over $k$. Writing
$$E \otimes k^{sep} \cong \underset{i \in S(\rho)}\oplus \underset{j =
1}{\overset{\rho(i)}{\oplus}} e_{i,j} k^{sep},$$ we have an action of
$G$ on the idempotents $e_{i,j}$. One may check that the idempotents
$e_{i,j}$ are permuted by $G$, and there is a correspondence between
the orbits of this action and the idempotents of $E$. In particular we have
\begin{Lem}
In the notation above, if $E$ is a subfield of $A$, then $|S(\rho)| =
1$.
\end{Lem}
\begin{proof}
$E$ is a field if and only if $G$ acts transitively on the set of
idempotents. On the other hand, this action must also preserve the
rank of an idempotent, which implies that all the idempotents have the
same rank.
\end{proof}

Therefore, if we are interested in studying the subfields of a central
simple algebra, we may restrict attention to partitions of the above
type.  If $m | n = deg(A)$, we write
$$\et_m(A) = \et_{[\frac{n}{m}, \frac{n}{m}, \ldots, \frac{n}{m}]}.$$

Note that every subfield of dimension $m$ is represented by a
$k$-point of $\et_m(A)$, and in the case that $A$ is a division
algebra, this gives a 1-1 correspondence. In particular, elements of
$\et_n(A)(F)$ are in natural bijection with the maximal subfields of
$A$. We will now show that this variety is rational and
R-trivial. This argument will be a geometric analog of one in
\cite{KraSa}.

If $a \in A$ is an element whose characteristic polynomial has
distinct roots, then the field $k(a)$ is a maximal subfield of $A$.

\begin{Thm} \label{unirational}
Let $U \subset A$ be the Zariski open subset of elements of $A$ whose
characteristic polynomials have distinct roots. Then there is a
dominant rational map $U \ra \et_n(A)$ which is surjective on
$F$-points.
\end{Thm}
\begin{proof}
Let $\til{U}$ be the degree $n!$ \'etale cover of $U$ whose
$k^{sep}$-points consists of pairs $(r, a)$, where $a \in U$, and $r =
(r_1, \ldots r_n)$ is an ordered $n$-tuple of distinct roots of the
chacteristic polynomial $\xi(a)$. One may check that the
indecomposable idempotents in the \'etale extension $k(a)$ are given
by the elements:
$$e_i = \underset{j \neq i}{\prod} \ \frac{a - r_j}{r_i - r_j}.$$ We
may then define a morphism $\til{U} \ra V_{\frac{n}{m}}(A)^d$ by
taking $(a, r)$ to $(e_1 A, \ldots e_n A)$. We compose this with the
quotient map and obtain a morphism $\til{U} \ra \et_n(A)$. Since this
morphism is constant on the fibers of $\til{U} \ra U$, this map
descends to a map $U \ra \et_n(A)$.

By the description above, this is the morphism which takes an element
of $A$ to the subfield which it generates. Since every \'etale
subfield of $A$ can be generated by a single element, this morphism is
surjective on $F$-points. Since this also holds after fibering with
the algebraic closure, it follows also that this morphism is
surjective at the algebraic closure and hence dominant.
\end{proof}

\begin{Cor}
Suppose $k$ is an infinite field. Then $\et_n(A)$ is a rational variety.
\end{Cor}
\begin{proof}
Suppose $L \subset A$ is a linear affine subvariety of $A$. Then from
remark \ref{grass} the $n$-dimensional \'etale subfields of $\et_n(A)$
which intersect $L$ transversely form an open subvariety $V$ of
$\et_n(A)$.

Choose a maximal \'etale subalgebra $E \subset A$, and let $L$ be a
complementary subspace of dimension $n^2 - n = dim(A) - dim(E)$.
Choose $a \in E$ with distinct eigenvalues generating $E$.  Note that
by construction, $(L + a) \cap E = a$.

Let $U \subset L + a$ be the dense open set of elements with distinct
characteristic roots, and $\psi$ be the restriction of the morphism
from theorem \ref{unirational} to $U$. I claim that this is a
birational isomorphism from $U$ to $\et_n(A)$. To see this note that
for $E' \in V$, there is a unique element $b \in E' \cap (L + a)$, if
we define $\phi(E') = b$, we obtain a morphism $V \overset{\phi}{\ra}
U$, which is a birational inverse to $\psi$.
\end{proof}

\begin{Lem}
Let $A$ be a degree $4$ central simple $k$-algebra. Then $V_2(A)$ is
isomorphic to an involution variety $V(B, \sigma)$ of a degree $6$
algebra with orthogonal involution $\sigma$.
\end{Lem}
\begin{proof}
Consider the map 
$$Gr(2,4) \ra \P^5$$, given by the Pl\"uker embedding. Fixing $V$ a
$4$ dimensional vector space, we may consider this as the map which
takes a $2$ dimensional subspace $W \subset V$ to the $1$ dimensional
subspace $\wedge^2W \subset \wedge^2V$. This morphism gives an
isomorphism of $Gr(2,4)$ with a quadric hypersurface. This quadric
hypersurface may be thought of as the quadric associated to the
bilinear form on $\wedge^2V$ defined by $<\omega_1, \omega_2> =
\omega_1 \wedge \omega_2 \in \wedge^4 V \cong F$. Note that one must
choose an isomorphism $\wedge^4 V \cong F$ to obtain a bilinear form,
and so it is only defined up to similarity. Nevertheless, the quadric
hypersurface and associated adjoint (orthogonal) involution depend
only on the similarity class and are hence canonically defined.

Since the Pl\"uker embedding defined above is clearly $PGL(V)$
invariant, using \cite{Ar:BS}, for any degree $4$ algebra $A$ given by a
cocycle $\alpha \in H^1(k, PGL_4)$, we obtain a morphism:
$$V_2(A) \ra V(B),$$ where $B$ is given by composition of $\alpha$
with the standard representation $PGL(V) \ra PGL(V \wedge V)$. By
\cite{Ar:BS} this implies that $B$ is similar to $A^{\otimes2}$ in
$Br(k)$. Also, it is easy to see that the quadric hypersurface and
hence the involoution is $PGL_4$ invariant, and hence descends to an
involution $\sigma$ on $B$. We therefore obtain an isomorphism
$V_2(A) \cong V(B, \sigma)$ as claimed.
\end{proof}

\begin{Lem}
Suppose $B$ is a central simple $k$ algebra with orthogonal involution
$\sigma$. Then there is a natural birational morphism $V_2(B) \ra
S^2V(B, \sigma)$ which is surjective on $k$ points. In particular,
$S^2V(B, \sigma)$ is a rational variety.
\end{Lem}
\begin{proof}
The proof is along similar lines to the previous proof. Let $V$ be a
vector space with a quadratic form $q$, and let $Q \subset \P(V)$ be
the associated quadric hypersurface. Consider the map
$$Gr(2,V) \ra S^2Q,$$ defined by taking a projective line in $\P(V)$
to its intersection points with $Q$. Note that this morphism is a
birational morphism. If we let $O_q(V)$ be the orthogonal group for
the quadratic form $q$, it is easy to see that the map described above
is $O_q(V)$ invariant.

To say that an algebra $B$ has an orthogonal involution is to say that
we have chosen a representing cocycle $\alpha \in H^1(k,
PO_q(V))$. Therefore, the above map will be invariant under any Galois
action defined by $\alpha$, and we will get a morphism
$$V_2(B) \ra S^2 V(B, \sigma).$$
\end{proof}

Putting these lemmas together yeilds the following corollary:

\begin{Cor}
Suppose $A$ is a degree $4$ central simple $k$ algebra. Then
$\et_2(A)$ is a rational variety.
\end{Cor}

\bibliographystyle{alpha}

\bibliography{citations}

\end{document}